\documentclass[a4paper, 11pt]{article}

\usepackage{amssymb}
\usepackage{amsmath}
\usepackage{graphicx, color}
\usepackage{pdfsync}

\def\z{\frak z}

\def\c{{\cal C}_z}

\def\nn{\nonumber}

\def\cut{\operatorname{cut}}

\def\tr{\operatorname{tr}}

\def\flecha{\longrightarrow}

\def\z{\frak z}

\def\ds{\displaystyle}
\def\fracc{\displaystyle \frac}

\def\s{{\rm s}_\lambda}
\def\c{{\rm c}_\lambda}
\def\ta{{\rm ta}_\lambda}

\def\cc{{\rm c}_{4\lambda}}

\def\ema{M^{n-q}_\lambda}

\def\emaz'd{M^{n-q-d}_{\lambda \frak z'}}

\def\emafz'd{M^{n-q-d}_{\lambda f(z')}}

\def\ema0{M^{n-1}_{\lambda 0}}
\def\cepe{\mathbb C P^{n}(\lambda)}
\def\cepeu{{\mathbb C P^{1}(\lambda)}}
\def\cepeq{{\mathbb C P^{q}(\lambda)}}
\def\qq{{Q^{n-1}(\lambda)}}

\def\pp{\frak P}

\def\la{\sqrt\lambda}
\def\ce{\mathbb C}

\def\<{\langle}
\def\>{\rangle}
\def\({\left(}
\def\){\right)}

\def\re{\mathbb R}
\def\${\sqrt\lambda}
\def\ce{\mathbb C}

 at 20truept
 at 24truept
 at 14truept
 at 16truept
\newtheorem{teor}{Theorem}

\newtheorem{nota}[teor]{Remark}

\numberwithin{teor}{section}

\begin{document}

\title{Bounding the first Dirichlet eigenvalue of a tube around a complex submanifold of $\cepe$ by the degrees of the polynomials defining it}
\author{M. Carmen Domingo-Juan and Vicente Miquel}
\date{}
\maketitle{}

\begin{abstract} We obtain  upper bounds for the first Dirichlet eigenvalue of a tube around a complex submanifold $P$ of $\cepe$ which depends only on the radius of the tube, the degrees of the polynomials defining $P$ and the first eigenvalue of some model centers of the tube. The bounds are sharp on these models. Moreover, when the models used are $\cepeq$ or $\qq$ these bounds also give gap phenomena and comparison results.
\end{abstract}

\section{Introduction}

This paper concerns with classical Dirichlet eigenvalue problem
\begin{equation}\label{DP}
\Delta f = \mu f  \ \text{ on } \ M \qquad \text{ and } \qquad 
f= 0 \ \text{ on }\ \partial M
\end{equation}
on a  connected compact Riemannian manifold $M$ with boundary $\partial M$. An active area of research in this problem is to determine \lq\lq interesting'' upper or lower bounds for the first eigenvalue $\mu_1(M)$ corresponding to equation \eqref{DP}. Here  \lq\lq interesting'' means sharp or related with some special properties of the space $M$. As examples of the work in this direction we have the results of S. Y. Cheng (\cite{Cg}), R. Reilly (\cite{Re}), M.Gage (\cite{Ga}), A. Kasue (\cite{Ks}), J. M. Lee(\cite{Lee}), F. Gim\'enez, A. Lluch, V. Palmer and the second author (\cite{GiMi, LlM, MP1, MP2}), G. P. Bessa and J. F. Montenegro (\cite{BM}). In these quoted papers, $M$ is a tube around some compact submanifold of some Riemannian manifold or a manifold with boundary with special bounds on the curvatures of $M$ and $\partial M$. The bounds obtained for $\mu_1(M)$ are usually the corresponding to some special model tube and, many times, the euqality characterizes this model tube. 

Very close to the the problem of obtaining bounds for $\mu_1(M)$ is the obtention of estimates for the volume of $M$. The deep relation between these problems is well known and some of the papers quoted before show or explore it. Concerning the volume, A. Gray (see \cite{Gra1,Gra2,Gr}) has shown that the volume of a tube around a complex submanifold of $\cepe$ can be stated in terms of the degrees of the polynomials defining the center of the tube. It is natural, then, to ask for some relation between the first Dirichlet eigenvalue of the tube and the degrees of the polynomials defining its center. Another hope for getting such kind of relation are the recent results by Colbois, Dryden and El Soufi (\cite{CDES}) where they obtain bounds of the first eigenvalue for the closed problem on algebraic submanifolds in the euclidean space in terms of the degrees of the polynomials defining them. The problem was addressed in \cite{DoLlMi}, where A. Lluch and the authors obtained a bound of the first Dirichlet eigenvalue of a tube around a complex curve in the complex projective space in terms of the degrees of the polynomial defining the complex curve. Our aim in this paper is to complete the work started in \cite{DoLlMi} by studying tubes around complex submanifolds of any dimension and , also, giving bounds of different flavor obatining the bound by comparison with the eigenvalues of tubes around all the possible homogeneous complex submanifolds for which the eigenfunction corresponding to $\mu_1(M)$ are radial.

We state now our results with more precision. Let $\cepe$ be the complex projective space of holomorphic sectional curvature $4 \lambda$, and let  $\wp:\ce^{n+1} \flecha \cepe$ be the canonical projection. Chow's Theorem states that every complete complex submanifold $P^q$ of $\cepe$ (of complex dimension $q$)   is the image by $\wp$ of the set of zeroes of $n-q$ homogeneous polynomials  of degrees $a_{q+1}, ..., a_n$. For such a submanifold, we shall denote by  $P_{\rho}$ the tube of radius $\rho > 0$ around $P$  and by $\partial P^q_{\rho}$ its boundary.   We shall always consider $\rho$ lower than the cut distance $\cut(P)$ from $P^q$. 
 In \cite{DoLlMi} it is studied the Dirichlet eigenvalue problem  \eqref{DP} for $M= P^1_\rho$, for $P$ a complex curve of $\cepe$ and it is obtained an upper bound of $\mu_1(P_\rho)$ of the form 
 \begin{align}
\mu_1(P^1_{\rho}) &\le   \mu_1(\cepeu_{\rho}) - M^1(\rho, n, a_2, \dots , a_q),  \label{DLM1}
\end{align} 
where $\cepeu$ is embedded as a totally geodesic complex submanifold of $\cepe$,  \\ $M^1(\rho, n, 1, \dots , 1)=0$  and $M^1(\rho, n, a_2, \dots, a_q)>0$ for $(a_2, \dots , a_n) \not= (1, \dots , 1)$. Moreover, the equality is attained if and only if $P^1=\cepeu$ (that is, $(a_2, \dots , a_q)=(1, \dots , 1)$).

This bound address the problem of relating $\mu_1(P_r)$ with the degrees of the polynomials defining $P$. Moreover it is also a comparison theorem with $\mu_1(\cepeu_r)$ and shows a gap phenomenon for $\mu_1(P^1_{\rho})$ between the case $P^1=\cepeu$ (which corresponds to $a_2 = \dots = a_n=1$) and the other complex submanifolds $P$, and states that the gap is measured by the degrees $a_s$ (a similar gap phenomenon occcurs in the study of the closed eigenvalue problem for complex submanifolds by J. P. Bourgugnon, P. Li and S. T. Yau (\cite {BLY})).

In this paper we address the problem of getting bounds of the same nature that \eqref{DLM1} for higher dimensions of $P^q$ ($q\ge1$) and  comparing with  tubes around some model complex submanifolds. These models, denoted by $\pp$, are the complex submanifolds of $\cepe$ having constant normal curvatures (constant means here that it does not depend on the point nor on the direction). They were classified by Kimura (\cite{Ki}) and are listed at the end of section 2.  As in \cite{DoLlMi}, we use the deep ideas in the work of A. Gray (\cite{Gra1}, \cite{Gra2}, \cite{Gr}) to get theorems. But here we use them more form the root, which gives simpler computations and also more general results. Of course, another ingredient is the work of Kimura (\cite{Ki}). We shall prove:

\begin{teor}\label{MT1} For $1\le q \le n-1$, $q=\dim_\ce P = \dim_\ce\pp$, the first eigenvalue $\mu_1(P^q_\rho)$ of the Dirichlet eigenvalue problem \eqref{DP} satisfies the inequalities
\begin{align}\label{mu1}
\mu_1(P^q_{\rho}) &\le   \mu_1(\pp_{\rho}) + M_\pp(\rho, n, q, a_{q+1}, ..., a_n),  \end{align} 
where $M_\pp(\rho, n, q, a_{q+1}, ..., a_n)$  are  well defined constants which depends only on $\pp, \rho, n, q,$ $a_{q+1}, ..., a_n$. 
Moreover: 
\begin{enumerate}
\item $M_\cepeq\le 0$ and the equality holds if and only if $P=\cepeq$, which is equivalent to $\sum_{i=q+1}^n a_i =n-q$, then \eqref{mu1} also gives  a gap between $\mu_1(\cepe_\rho)$ and the corresponding eigenvalues of tubes with the same radius around complex submanifolds defined by polynomials of higher degree.
\item  When $q = n - 1$ and $\pp = Q^{n - 1}(\lambda)$ (the complex hyperquadric), there is a $\rho_0$, $0<\rho_0\le \cut(P)$, depending on the degree $a_n$ of the polynomial defining $P$,  such that, for every $\rho < \rho_0$, 
$M_{Q^{n - 1}(\lambda)} \leq 0$ if $P \ne \cepe$ (which equivalent to $a_n > 1$) and, in this case, equality holds if and only if the polynomial defining $P$ has degree $2$.  As a consequence, for every $\rho < \rho_0$:
 \begin{enumerate}
\item $Q^{n-1}(\lambda)$ gives the biggest first eigenvalue of problem \eqref{DP} among all the tubes around a complex hypersurface defined by an homogeneous polynomial of degree $2$, and 
\item For all the complex hypersurfaces defined by polynomials of degree $\ge 3$, there is a gap between $\mu_1(\qq_{\rho})$ and  $\mu_1(P_\rho) $.
\end{enumerate}
 \end{enumerate}

\end{teor}

\begin{nota}
We think that $\rho_0= \cut(P)$, but we have not enough precise bounds for $\rho_0$ to assure it. For the other models $\pp$ different from $\cepe$ and $Q^{n-1}(\lambda)$ we have no control on the term $M_\pp$. This is what make us unable to obtain any kind of comparison theorem from the bound ($ \mu_1$) for these other models. 

If $\pp$ is not any of the models we have considered, the $\mu_1$-eigenfunction is not radial, then our method cannot give an  upper bound of $\mu_1(P_\rho)- \mu_1(\pp_\rho)$ depending only on , $\pp$, $n$, $q$,   $\rho$ and the degrees of the polynomilas defining $P$.
\end{nota}

{\bf Acknowledgments:} 
Second author has been partially supported  by DGI (Spain) and FEDER Project   MTM2010-15444 and the Generalitat Valenciana Project GVPrometeo 2009/099.

\section{Preliminaries on complex submanifolds and the tubes around them}

Given a complex submanifold $P$ of $\cepe$ of real dimension $2q$, we shall denote by   $r$ the distance to $P$ in $\cepe$. Let us denote by $\mathcal NP$ the normal bundle of $P$, by $A_\xi$ the Weingarten map of $P$ in the direction of $\xi\in \mathcal N P$, $|\xi|=1$. 
Moreover, we shall use the notations $\s$ and $\c$ for
$$\s(t) = \fracc{\sin(\sqrt{\lambda} t)}{\sqrt{\lambda}}, \qquad \c(t) = \cos(\sqrt{\lambda} t), \qquad \ta(t)=\frac{\s(t)}{\c(t)}$$
which satisfy the computational rules $\s'=\c$ and $\c^2 + \lambda \s^2=1$.

Since $P$ is a complex submanifold, given $\xi\in \mathcal N P$, $|\xi|=1$,  the Weingarten map $A_\xi$, has eigenvalues $k_1(\xi), -k_1(\xi), ..., k_q(\xi), -k_q(\xi)$.
The trace of  the Weingarten map $S(t$) of $\partial P_t$  is (cf. \cite{Gr}, page 125, formula (7.25)):
\begin{align}\label{trS}
&\tr S = 2 \s(\rho) \c(\rho)\ h_P(\rho)  - (2 n - 2q - 1) \frac{\c(\rho)}{\s(\rho)} 
 + \lambda \frac{ \s(\rho)}{c_{ \lambda}(\rho)},\\
 & \text{ where } h_P(\rho) =\sum_{i = 1}^q 
\frac{( \lambda + k_i^2)}{\c^2(\rho) - k_i^2  \s^2(\rho)} \label{hp}
\end{align}

On the other hand, if $f:\re \flecha \re$ is a $C^\infty$ function, one has (cf. \cite{MP1} for instance): 
\begin{equation}
 \label{radial}
 \Delta (f\circ r) = - f''\circ r + \tr S \ f'\circ r
\end{equation}
From now on we shall omit the writing of \lq\lq$\circ \ r$'' when it can be understood by the context. 
 
 \medskip
 
 The volume element $\omega$ of a tube $P_\rho$ in Fermi coordinates around $P$ can be written (cf. \cite{Gr}, page 125 formula (7.26)) as
 \begin{align}\label{ve}
 & \omega = \theta(p,\xi, r) \ d\xi\ dp\ dr, \text{ with } \theta(p,\xi, r) = \s^{2n-2q-1} \  \c \ v(p,\xi, r) , \\
 &\quad \text{where } v(p,\xi,r) =\prod_{j=1}^q (\c^2- \s^2 \ k_j(\xi)^2) =   \c^{2 q } \prod_{j=1}^q (1- \ta^2 \ k_j(\xi)^2) , \nn
 \end{align}
 where $dp$ and $d\xi$ denote, respectively, the volume elements of  $P$ and $S^{2n-2q-1}$. We remark that, for $p$ and $\xi$ fixed, the first positive value of $r$ where $v(p,\xi,r)$ (then $\theta(p,\xi,r)$) vanishes is lower than $\cut(P)$.
 
Developing the product or the determinant in the above formula, one obtains 
\begin{align}
&\prod_{i=1}^q (1- \ta^2 \ k_i(\xi)^2) = \sum_{c = 0}^q \Psi _{2 c} (\xi, \dots , \xi)\  \ta^{2c}   , \text{ with }\label{scaner2}\\
& \quad \Psi _{2c} (\xi, \dots , \xi) = (-1)^c \sum_{
\begin{array}{c}
i_1, \dots , i_c = 1 \\
 i_1 < \dots < i_c
\end{array}
}^q k_{i_1}^2(\xi) \dots k_{i_c}^2(\xi) 
\label{scaner4}
\end{align}
And the $\psi_{2c}$ satisfy (cf. \cite{Gr}, pages 65 and 125):
\begin{align}
\int_{S^{2n-2q-1}} \Psi _{2c} (\xi, \dots , \xi) \ d\xi =: &I_{2c} (\Psi_{2c} )  = a(c) \ \mathcal C ^{2c} ((R^P - R^{\cepe})^{c}), \label{scaner9}
\end{align}
where 
\begin{align}\label{ac}
a(c) = \ds{\frac{2 \pi ^{n - q}}{c! \ (2c)! \ 2^{c} \  (n - q + c-1)!}} 
\end{align}
 and $\mathcal C ^{2c} ((R^P - R^{\cepe})^{c})$ is a contraction of the curvature operator $(R^P - R^{\cepe})^{c} $ which is related with the  Chern form  $\overline{\gamma}_c= \gamma _c (R^P - R^{\cepe})$ of the curvature operator  $(R^P - R^{\cepe})$ by 
\begin{align}
 b(c) := \ds  \int_{P} \ \mathcal C ^{2c} ((R^P - R^{\cepe})^{c})  \ dp = \frac{i! \ (2c)! \ (2 \pi )^{c} }{(q - c)!} \int_{P} \ \overline{\gamma}_c \wedge F^{q - c}  \ dp. \label{mipg15}
\end{align}
One can look at \cite{Gr} page 56 for a precise definition of $\mathcal C ^{2c} $ and page 88 for the definition of $\gamma_c$. On the other hand, $\overline \gamma_c$ is related with the degrees $a_j$ of the polynomials defining $P$ (cf. \cite{Gr} page 141) by 
 \begin{align}
  \left[ \overline{\gamma} \right] = 
\left[ \gamma (R^P - R^{\cepe})  \right] = \left[ \frac{1}{\ds \prod_{j = q+1}^n\left( 1 + \frac{(a_j - 1) \lambda }{\pi } F  \right) } \right],
\end{align} 
where $[ \cdot{}  ]$ denotes the cohomology class of the corresponding differential form and the $\overline \gamma_c$ are defined from $\overline \gamma$ by 
\begin{align*}
\overline{\gamma} 
= 1 + \overline{\gamma}_1  + \dots + \overline{\gamma}_q  + \dots \ . 
\end{align*}
From which it follows that
\begin{align*}
& 1 = \left( 1 + \overline{\gamma}_1 + \dots + \overline{\gamma}_q  + \dots \right) 
\prod_{j = q+1}^{n} \(1 + (a_j - 1)  \frac{\lambda}{\pi} F\) ,
\end{align*}
then, for $c= 1, \dots , q$ ,
\begin{align}
  \overline{\gamma}_c  =  (- 1)^{c}  \beta_c  \( \frac{\lambda}{\pi} F\)^c \quad \text {where } \quad \beta_c = \sum_{\begin{array}{c}
j_1, \dots , j_{c} = q+1 \\
 j_1 \leq \dots \leq j_c 
\end{array}}^{n} (a_{j_1} - 1) \dots (a_{j_c} - 1) 
 \label{mipg17-Ib}
\end{align}
By substitution of   \eqref{mipg17-Ib} in \eqref{mipg15} one obtains 
\begin{align}
b(c) & = (- 1)^{c}  \frac{c! \ (2c)! \ (2 \pi )^{c} \  \lambda ^{c} \   \beta_c  }{(q - c)! \ \pi ^{c} } \int_P 
F \wedge \overset{\underset{\smile}{q}}{\dots} \wedge F \ dp 
= (- 1)^{c}  \frac{c! \ (2c)! \ 2^{c} \ \ q! }{(q - c)! } \lambda ^{c} \   \beta_c \  vol(P)  .
 \label{mipg18}
\end{align}

In \cite{Ki}, Kimura classified all the complex submanifolds of $\cepe$ whose principal curvatures are constant in the sense that they depend neither on the point of the submanifold nor on the normal vector.  They are: 
\begin{align}\label{trla} &\bullet \text{ Totally geodesic $\cepeq$. It has $k_i=0$}\\
\label{trQ} &\bullet\text{ The complex hyperquadric $Q^{n-1}(\lambda)$, where $k_i(\xi) = \sqrt{\lambda}$ for $i=1, ..., n-1$.}\\
&\bullet \label{trC}  \text{ $\mathbb C P^{1}(\lambda) \times \mathbb C P^{m - 1}(\lambda) \subset 
\mathbb C P^{2m - 1}(\lambda)$ for $m \geq 3$ (then $n= 2m-1\ge 5$}\\
&\qquad \text{ and $q=m$),  where $k_1 = k_2= \la$, $k_3 = ... = k_m =0 $}\nn\\
&\bullet\label{trU}  \text{ $SU(5) / S(U(3) \times U(2)) \subset \mathbb C P^{9}(\lambda)$ (then $n=9$, $q=6$), }\\
&\qquad \text{where $k_1 = k_2=k_3 =k_4= \la$, $k_5 =  k_6 =0 $. }\nn\\
&\bullet \label{trO} \text{ $SO(10) / U(5) \subset \mathbb C P^{15}(\lambda)$ (then $n=15, q=10$),}\\
&\qquad \text{  where $k_1 = ... =k_6= \la$, $k_7 = ... =  k_{10} =0 $}\nn
 \end{align}

 \section{Proof of Theorem \ref{MT1}}
 
 We shall denote by $S_\pp$ the Wingarten maps of the tubular hypersurfaces centred at a model comparison $\pp$ and by $f_\pp$ an eigenfunction corresponding to $\mu_1(\pp_\rho)$. According to \eqref{radial}, $f_\pp$ is the solution of the equation
\begin{align}
 \label{CPq}
 \begin{array}{l}
- f_\pp'' + \tr S_\pp\ f_\pp' =  \mu_1(\pp_\rho)  f_\pp ,  \qquad 
f_\pp(\rho) = 0. \qquad 
f_\pp' (0) = 0. 
\end{array}
\end{align}
 This function satisfies the inequalities (cf. \cite{LlM}  or \cite{MP1} for instance)
 \begin{equation}\label{signff}
 f_\pp >0 \text{ on } [0,\rho[ \quad \text{ and } \quad  f_\pp' <0 \text{ on } ]0,\rho].
 \end{equation}
   
  We shall apply the Raileigh's theorem 
using   $f_\pp\circ r$  as a test function, then  
 \begin{equation}
 \label{cociente}
\mu _1 (P_{\rho})  \le  \frac{\int _{P_{\rho}} f_\pp(\Delta f_\pp) }{\int _{P_{\rho}} f_\pp^2}
 \end{equation}
Let us compute the right hand side of the above inequality. From \eqref{trS},
\eqref{radial} and \eqref{CPq}, we get
 \begin{align}
 \Delta (f_\pp \circ r) & = -   f_\pp''\circ r  + \tr S_\pp\ f_\pp' \circ r+ \( \tr S\ - \tr S_\pp \) f_\pp' \circ r \nn\\
 &= \mu_{1}(\pp_\rho) f_\pp+ 2 \s\circ r\ \c\circ r \( h_P - h_\pp \) f_\pp'   \label{Delta}
\end{align}
 
From \eqref{cociente}, \eqref{Delta}, \eqref{trla} and \eqref{trQ} one gets
\begin{align}
& \mu_1(P_{\rho}) \leq   \mu_{1}(\pp_\rho)
 +  \frac{\ds\int _{P_{\rho}} f_\pp\ 2\ \s \ \c  
 \( h_P - h_\pp\)   f_\pp' \ \omega}{\ds\int _{P_{\rho}} f_{\pp}^2\ \omega} \nonumber \\ 
& = \mu_{1}(\pp_\rho) 
 +  \frac{\ds\int _0^{\rho} \int_{P} \int_{S^{2n-2q-1}}   2 f_{\pp} f'_{\pp}\  \s^{2n-2q} \ \c^2 \ 
 \( h_P - h_\pp\)	 \   \ v(p,\xi, r)\ d\xi\ dp\ dr}
{\ds\int _0^{\rho}    \int_{P} \int_{S^{2n-2q-1}} f^2_{\pp}\ \theta(p,\xi, r) \ d\xi\ dp\ dr}, \label{pg11}
\end{align}
where we have used the expression \eqref{ve} of the volume element $\omega$ of $P^q_\rho$ in Fermi coordinates.

Let us work with the integrand of the numerator in \eqref{pg11}. Using \eqref{ve}, \eqref{scaner2} and \eqref{scaner4} we get
\begin{align}\label{intnu}
h_P\ v(p,\xi,r) &= \left( \sum_{i = 1}^q \ \frac{ k_i(\xi)^2 + \lambda}{\c^2 - k_i(\xi )^2  \s^2}\right) 
	\prod_{j=1}^q(\c^2- k_j(\xi)^2 \s^2)  \nn \\
&=  \sum_{i = 1}^q (k_i(\xi )^2 + \lambda) \ \c ^{2q-2} \prod_{\begin{array}{c}
j = 1 \\
j \not= i 
\end{array}}^q (1 - k_j(\xi )^2  \ta^2) 
\end{align}
\begin{align}\label{intnupp}
h_\pp\ v(p,\xi,r) &= h_\pp	\prod_{j=1}^q(\c^2- k_j(\xi)^2 \s^2)  
=  h_\pp  \ \c ^{2q} \sum_{c=0}^q  (-1)^c\ \ta^{2c}\sum_{
\begin{array}{c}
i_1, \dots , i_c = 1 \\
 i_1 < \dots < i_c
\end{array}
}^q k_{i_1}^2 \dots k_{i_c}^2 
\end{align}

But direct computation gives
\begin{align}
& \sum_{i = 1}^q k_i(\xi )^2 \prod_{\begin{array}{c}
j = 1 \\
j \not= i 
\end{array}}^q (1 - k_j(\xi )^2  \ta^2) = \sum_{c = 1}^q (-1)^{c - 1} \ c \  \ta^{2c - 2} \sum_{
\begin{array}{c}
i_1, \dots , i_c = 1 \\
 i_1 < \dots < i_c
\end{array}
}^q k_{i_1}^2 \dots k_{i_c}^2 
\label{scaner1}
\end{align}
and
\begin{align}
& \sum_{i = 1}^q  \lambda \prod_{\begin{array}{c}
j = 1 \\
j \not= i 
\end{array}}^q (1 - k_j(\xi )^2  \ta^2) =  \lambda  \sum_{c = 1}^q (-1)^{c - 1} \ (q-c+1) \  \ta^{2c - 2} \sum_{
\begin{array}{c}
i_2, \dots , i_c = 1 \\
 i_2 < \dots < i_c
\end{array}
}^q k_{i_2}^2 \dots k_{i_c}^2 
\label{scaner1a}
\end{align}

From \eqref{intnu}, \eqref{scaner1},  \eqref{scaner1a}, \eqref{scaner4}, \eqref{scaner9} and \eqref{mipg15}   it follows
\begin{align}\label{numsrp}
\int_{P} &\int_{S^{2n-2q-1}}   h_P \ v(p,\xi, r)  \nn \\
	& = \c ^{2q-2} \( -  \sum_{i=1}^q i\ \ta^{2i-2} \ a(i) \ b(i) +  \lambda\   \sum_{i=1}^q (q-i+1) \ \ta^{2i-2} \ a(i-1) \ b(i-1) \)
\end{align} 
and
\begin{align}\label{numsrfp}
\int_{P} &\int_{S^{2n-2q-1}}   h_\pp \ v(p,\xi, r)   =   h_\pp \ \c ^{2q} \sum_{i=0}^q  \ta^{2i} \ a(i) \ b(i)
\end{align} 

Now, after the substitution of \eqref{numsrp} and  \eqref{numsrfp} in the numerator of  \eqref{pg11}, having into account that $2 f_{\pp} f'_{\pp}  = ( f_{\pp}^2 )'$, we compute for that numerator:

\begin{align}
&\ds\int _0^{\rho}  \int_{P} \int_{S^{2n-2q-1}}   2 f_{\pp} f'_{\pp}\  \s^{2n-2q} \ \c^2 \ 
 \( h_P - h_\pp\)	 \   \ v(p,\xi, r)\ d\xi\ dp\ dr \nn \\
&=   \int _0^{\rho}   ( f_{\pp}^2 )'  \s^{2n-2q} \ \c^{2q} \  
\(  -  \sum_{i=1}^q i\ \ta^{2i-2} \ a(i) \ b(i) +  \lambda\   \sum_{i=1}^q (q-i+1) \ \ta^{2i-2} \ a(i-1) \ b(i-1)   \right.\nn \\
& \qquad \qquad \qquad \qquad \qquad \qquad \qquad \qquad\qquad \left. -  h_\pp \ \c ^{2} \sum_{i=0}^q  \ta^{2i} \ a(i) \ b(i)  \) \ dr \nonumber \\ 
&=   \int _0^{\rho}   ( f_{\pp}^2 )' \  \s^{2n-2q} \ \c^{2q} \  
  \sum_{i=0}^q \(  - i\ \ta^{2i-2} +  \lambda\   (q-i) \ \ta^{2i}  -  h_\pp \ \c ^{2}   \ta^{2i} \) a(i) \ b(i)   \ dr =:\mathcal N \label{numint} 
  \end{align}
  But in all our comparison models $\pp$ the possible values for $k_i$ are $\pm\la$ or $0$. Let $\z$ be the number of normal curvatures with value $\la$. With this conventon, it follows from \eqref{hp} that the expression for $h_\pp\ \c^2$ is 
  \begin{equation}
h_\pp\ \c^2=   \z\ \lambda \frac{1}{\cc} + q \lambda
  \end{equation}
  We continue with the above computation \eqref{numint}
  \begin{align}
 \mathcal N &=   \int _0^{\rho}   ( f_{\pp}^2 )'  \s^{2n-2q} \ \c^{2q} \  
  \sum_{i=0}^q \(  - i\ \ta^{2i-2} \c^{-2}+  \lambda\   q \ \ta^{2i}  -  \( \z\ \lambda \frac{1}{\cc} + q \lambda\)   \ta^{2i} \) a(i) \ b(i)   \ dr \nonumber \\ 
&  =   \sum_{i=0}^q  \int _0^{\rho}   ( f_{\pp}^2 )'   \  
  \( - i\  - \z \lambda \fracc{\s^2}{\cc} \) \s^{2(n-q+i-1)} \c^{2(q-i)}\ a(i) \ b(i)   \ dr \nonumber \\ 
 & =   \sum_{i = 0}^q  \int_0^{\rho} ( f_{\pp}^2 )'   \  
  \( \fracc{ - i\ \c^2  + (-\z + i) \lambda \s^2}{\cc} \) \s^{2(n-q+i-1)} \c^{2(q-i)}\ a(i) \ b(i)   \ dr \nn \\
  & =   \int_0^{\rho} ( f_{\pp}^2 )' \fracc{1}{\cc} \s^{2(n-q-1)} \c^{2 q }\sum_{i = 0}^q     \  
  \( - i\ \c^2  + (i - \z) \lambda \s^2 \) \ta^{2 i}\ a(i) \ b(i)   \ dr
 \label{mipg6}
\end{align}

On the other hand, for the denominator of \eqref{pg11}, from  \eqref{ve}, the remark after \eqref{ve}, \eqref{scaner2}, \eqref{scaner9} and \eqref{mipg15} it follows: 
\begin{align}\label{den}
0 < \int_0^\rho \int_{P} \int_{S^{2n-2q-1}}   & f_\pp^2	\theta(p,\xi, r)\ d\xi\ dp dr =  \int_0^\rho   f_\pp^2\ \s^{2n-2q-1} \c^{2q+1} \sum_{i=0}^q \ \ta^{2i} \ a(i) \ b(i) \ dr\nn \\
&= \sum_{i = 0}^q \ a(i) \ b(i)  \int_0^{\rho}  f^2_{\pp} \s^{2n - 2q  + 2i -1} \c^{2q - 2i + 1} \ dr 
\end{align} 
Now we shall use the notation 
\begin{equation}
\mu_\pp = f_\pp^2  \  \s^{2n - 2q -1} \ \c^{2q + 1}\ \text{ and } \  \nu_\pp = (f_\pp^2)'  \  \frac{\s^{2n - 2q -1} \ \c^{2q}}{\cc}.
\end{equation}

Substitution of \eqref{mipg6} and \eqref{den} in \eqref{pg11} gives 
\begin{align}\label{mipg7}
& \mu_1(P_{\rho}) \leq 
\mu_1(\pp_{\rho}) \\ 
& \qquad +  \frac{\ds  \int_0^{\rho}  \sum_{i = 0}^q     \  
  \( - i\ \c^2  + (i - \z) \lambda \s^2 \) \ta^{2 i}\ a(i) \ b(i) \nu_\pp  \ dr }
{\ds  \sum_{i = 0}^q \int_0^{\rho}   \s^{ 2i } \c^{- 2i } \ a(i) \ b(i)\ \mu_\pp\ dr },  \nn
\end{align}

From the formulae \eqref{ac} and \eqref{mipg18} for $a(i)$ and $b(i)$ it follows that 
\begin{align}
\fracc{a(i)\ b(i)}{a(0)\ b(0)} = \fracc
{\ds   \  (- 1)^{i}   \ (n - q -1)! \ q! \ \lambda ^{i} \   \beta_i  }{\ds   (q - i)! \  (n - q + i-1)!   } =  (-1)^i \fracc{\binom{n-1}{q-i}}{\binom{n-1}{q}} \lambda^i \ \beta_i. \label{forbi}
\end{align}
Then 
\begin{align}
&\frac{\binom{n-1}{q}}{ a(0) b(0)} \sum_{i = 0}^q     \  
  \( - i\ \c^2  + (i - \z) \lambda \s^2 \) \ta^{2 i}\ a(i) \ b(i) \nn \\
 &\qquad\quad=  \sum_{i = 0}^q     
  \( - i\ \c^2  + (i - \z) \lambda \s^2 \) \ta^{2 i}\ (-1)^i \binom{n-1}{q-i} \lambda^i \ \beta_i 
\nn \\
&\qquad\quad=  - \ \lambda\  \z\ \s^2  \binom{n-1}{q} \nn \\
&\qquad\qquad - \(-\s^2 + (1-\z) \lambda \s^2 \ta^2\) \lambda \binom{n-1}{q-1} \ \beta_1
\nn \\
&\qquad\qquad + \(- 2 \s^2 \ta^2 + (2-\z) \lambda \s^2 \ta^4\) \lambda^2 \binom{n-1}{q-2}  \ \beta_2
\nn\\
&\qquad\qquad - \(- 3 \s^2 \ta^4 + (3-\z) \lambda \s^2 \ta^6\) \lambda^3 \binom{n-1}{q-3}  \ \beta_3
\nn \\
&\qquad\qquad \dots \dots \dots \nn \\
&\qquad\qquad + (-1)^{q}\( - q \s^2 \ta^{2q-2} + (q-\z) \lambda \s^2 \ta^{2q}\) \lambda^q    \  \beta_q 
\nn\\
&\qquad\quad=  \sum_{i=0}^{q-1} (-1)^i \lambda^{i+1} \s^{2} \ta^{2i}\( (i-\z)\ \binom{n-1}{q-i} \beta_i +  (i+1) \binom{n-1}{q-i-1} \beta_{i+1} \) \nn\\
&\qquad\qquad \qquad +  (-1)^q (q-\z) \lambda \s^2 \ta^{2q} \lambda^q  \  \beta_q
\end{align}
Multiplyining numerator and denominator of \eqref{mipg7} by 
$\fracc{\binom{n-1}{q}}{a(0) \ b(0)} $, having into account \eqref{ac} and \eqref{mipg18}, we obtain

\begin{align} \label {mipg13b} 
& \mu_1(P_{\rho}) \leq 
\mu_1(\pp_{\rho})  \nn \\
&\qquad+   \frac{\ds \sum_{i = 0}^{q-1}  (-1)^i \lambda^{i+1} \((i-\z) \binom{n-1}{q-i} \beta_i + (i+1) \binom{n-1}{q-i-1} \beta_{i+1} \) B_i(\rho) +  (-1)^q (q-\z) \lambda \s^2 \ta^{2q} \lambda^q  \  \beta_q B_q(\rho)}
{ \ds \sum_{i = 0}^q  (-1)^i \binom{ n-1}{q-i} \lambda^i \beta_i  C_i(\rho )}  
\end{align}
where
 $B_i(\rho) = \ds \int_0^{\rho} \  \s^{2} \ta^{2i} \ \nu_\pp\ dr$ and $ \ds C_i(\rho) =  \int_0^{\rho} \  \s^{2i}  \c^{-2i}  \mu_\pp\ dr$.

When $\pp = \cepeq$, $\z=0$. To study the sign of  $M_{\cepeq}$ we consider its expression given by the second summand in \eqref{pg11}. In this case,  from the definition \eqref{hp} of $h_P$, one has $h_P-h_\cepeq = \sum_{i=1}^q \fracc{k_i^2}{\c^2 (\c^2 - k_i ^2 \s^2)} \ge 0$ for $r < \cut(P)$. Since also $\theta \ge 0$ for $r < \cut(P)$ and $f_{\pp}' <0$ on $]0,\rho]$, one has that $M_{\cepeq}\le 0$. On the other hand, the equality in \eqref{mu1} implies the equality in \eqref{cociente}, which implies that  $f_\cepeq$ is an eigenfunction with eigenvalue $\mu_1(\cepeq)$, which, from \eqref{Delta}, implies $h_P-h_\cepeq=0$, which, from the above expression only happen if $k_i=0$, that is, if $P=\cepeq$.

When $\pp=Q^{n-1}(\lambda)$, one has  $\z=q=n-1$ and there is only one polynomial defining $P$, with degree $a_n$, then $\beta_i = (a_n-1)^i$  and the numerator of $M_{Q^{n-1}(\lambda)}$ becomes
\begin{align}
\sum_{i = 0}^{n-2} & (-1)^i \lambda^{i+1} \((i-(n-1)) \binom{n-1}{n-1-i} (a_n-1)^i + (i+1) \binom{n-1}{n-1-i-1}(a_n-1)^{i+1} \) B_i(\rho) \nn \\
& = \sum_{i = 0}^{n-2}  (-1)^i \lambda^{i+1} \((i-(n-1)) \binom{n-1}{n-1-i} + (i+1) \binom{n-1}{n-1-i-1}(a_n-1) \) (a_n-1)^i B_i(\rho) \label{forty}
\end{align}
But 
\begin{align}
(i-(n-1))  \binom{n-1}{n-1-i} &= - (n-1-i) \binom{n-1}{n-1-i} \nn \\ &= - \fracc{(n-1) \cdots (n-i-1)}{i!} = - (i+1) \binom{n-1}{n-1-i-1},
\end{align}
which, substituted in  \eqref{forty} gives
\begin{align}
&\sum_{i = 0}^{n-2}  (-1)^{i} \lambda^{i+1} \fracc{(n-1) \cdots (n-i-1)}{i!} (-1+a_n-1)  (a_n-1)^i B_i(\rho)\nn \\
& = (a_n-2) \sum_{i = 0}^{n-2}  (-1)^{i} \lambda^{i+1} \fracc{(n-1) \cdots (n-i-1)}{i!}   (a_n-1)^i B_i(\rho)  \label{forty2}
\end{align}

Then $M_{Q^{n-1}(\lambda)}$ vanishes when $a_n=2$. As a consequence, among all the complex hypersurfaces $P$ defined by a polynomial of degree $2$, the complex hyperquadric gives the maximum value of $\mu_1(P)$. 

Now, let us study the sign of $M_{\qq}$. First, recall that its denominator is positive for $\rho < \cut(P)$ (as we noticed in \eqref{den}). To check the sign of the numerator, the observation that  when $2[(n-2)/2]+1 > n-2$,  $n-2 [(n-2)/2] -2 =0$, allows us to write the sum in \eqref{forty2} in the following way:
\begin{align}
\sum_{i = 0}^{n-2}  & (-1)^{i} \lambda^{i+1} \fracc{(n-1) \cdots (n-i-1)}{i!}   (a_n-1)^i B_i(\rho) \nn \\ 
&= \sum_{j=0}^{[(n-2)/2]} (-1)^{2j} \lambda^{2j+1} \fracc{(n-1) \cdots (n-2j-1)}{(2j)!}   (a_n-1)^{2j} \(B_{2j}(\rho) - \lambda \fracc{(n-2j-2)}{(2j+1)}   (a_n-1) B_{2j+1}(\rho)\)\nn \\
&= \sum_{j=0}^{[(n-2)/2]} \lambda^{2j+1} \fracc{(n-1) \cdots (n-2j-1)}{(2j)!}   (a_n-1)^{2j} \int_0^\rho \s^{2} \ta^{4j} \ \nu_\pp \( 1 - \lambda \fracc{(n-2j-2)}{(2j+1)}   (a_n-1)   \ta^{2} \) \ dr  \nn
\end{align}
which is negaitive for $\ds \rho \le \rho_1:= \min_{0\le j \le  [(n-2)/2]} \ta^{-1}\(\sqrt{\fracc{2j+1} {\lambda (n-2j-2)}}\)$, where $\ta^{-1}$ means the inverse function of $\ta$ with image in $[0,\pi/2\la[$. Then, taking $\rho_0 = \min\{\cut(P), \rho_1\}$, we have that $M_\qq <0$ for $a_n\ge3$. This gives a gap between $\mu_1(\qq_\rho)$ and $\mu_1(P_\rho)$ for all complex hypersurfaces defined by polynomials of degree $\ge 3$.

{\footnotesize

\bibliographystyle{alpha}

\begin{thebibliography}{99}
 
 \bibitem{BLY}  J. P. Bourguignon, P. Li and S. T. Yau,   
 Upper bound of the first eigenvalue for algebraic submanifolds, {\em Com. Math. Helv.} {\bf 69} (1994), 199-207
 
\bibitem{BM} G. P. Bessa and J. F. Montenegro, On Cheng's eigenvalue comparison theorem, {\em Math. Proc. Camb. Phil. Soc.}  {\bf 144} (2008), 673-682.

\bibitem {Cg} S.Y. Cheng,
Eigenvalue comparison theorems and its geometric applications,  {\em Math. Z.} {\bf 143} (1975),  289-297

\bibitem {CDES} B. Colbois, E. B. Dryden and A. El Soufi,
Bounding the eigenvalues of the Laplace-Beltrami operator on compact submanifolds,  {\em Bull. London Math. Soc.} {\bf 42} (2010),  96-108

\bibitem {DoLlMi} M. C. Domingo-Juan, A. Lluch and V. Miquel, 
Upper bounds for the first Dirichlet eigenvalue of a Tube around an algebraic complex curve of $\cepe$,  {\em Israel J. Math. } {\bf  183} (2011),   189--198.

\bibitem{Ga} M. Gage, Upper bounds for the first eigenvalue of the Laplace-Beltrami operator, {\em Indiana Univ. Math. J. } {\bf 29} (1980), 897--912.

\bibitem{GiMi}  F. Gim\'enez and V. Miquel, Bounds for the first
Dirichlet eigenvalue of domains in K\"ahler manifolds {\em Archiv der Math.}  {\bf 56} (1991)  370--375.

\bibitem{Gra1} A. Gray,  
 Volumes of tubes about K\"ahler submanifolds expressed  in terms of
Chern classes,{\em J. Math. Soc. Japan} {\bf 36} (1984) 23--35

\bibitem{Gra2} A. Gray,  Volumes of tubes about complex submanifolds of
complex projective space,  {\em Trans. Amer. Math.
Soc.} {\bf 291}, (1985) 437-449

\bibitem {Gr}  A. Gray {\em Tubes}, Second Edition,
Birkh\"auser,  Heidelberg, New York,
(2003).

\bibitem{Ks} A. Kasue, On a lower bound for the first eigenvalue of the Laplace operator on a Riemannian manifold, {\em Ann. Sci. Ëcole Norm. Sup. } {\bf 17} (1984), 31-44.

\bibitem{Ki}  M. Kimura, {Real Hypersurfaces and complex
submanifolds in complex projective space}. {\em Trans. Am. Math. Soc.} {\bf  296}
(1986)  137--149.

\bibitem{Lee} J.M. Lee, 
Eigenvalue comparison for tubular domains,  {\em Proc. Am. Math. Soc.} {\bf 109}, (1990) 843-848

\bibitem{LlM} A. Lluch, V. Miquel  Bounds for the first Dirichlet eigenvalue attained at an infinite
family of Riemannian manifolds 
{\em Geometriae Dedicata}, {\bf  61} (1996) 51--69

\bibitem{MP1} V. Miquel, V. Palmer, Mean curvature comparison for tubular hypersurfaces in
K\"ahler manifolds and some applications  {\em Compositio Math.} {\bf 86} (1993) 317--335.

\bibitem{MP2} V. Miquel, V. Palmer, Lower bounds for the mean curvature of hollow tubes around complex hypersurfaces and totally real submanifolds. {\em Illinois J. Math.} {\bf 39} (1995), 508--530.

\bibitem{Re} R.Reilly, Applications of the Hessian operator in a Riemannian manifold, {\em Indiana Univ. Math. J} {\bf 26} (1977), 459--472.

\end{thebibliography}

}

\vskip1truecm

{\small 

\begin{tabular}{ c c }
&Universidad de Valencia \\
&Departamento de Matem\'aticas para la Econom\'{\i}a y la Empresa \\
&Avda Tarongers s/n\\
& 46022-Valencia (Spain)\\
      \       & email: carmen.domingo@uv.es\\
      
& and \\

&Universidad de Valencia.\\
&Departamento de Geometr\'{\i}a y Topolog\'{\i}a\\
& Avda. Andr\'es Estell\'es, 1, 46100-Burjassot (Valencia) Spain \\
      \       & email: miquel@uv.es 
\end{tabular}

}

\end{document}